\documentclass[12pt,a4paper]{article}
\usepackage{amsthm}
\usepackage{amsmath}
\usepackage{amssymb}
\usepackage{amsfonts}
\usepackage{array}
\usepackage{hyperref,xcolor}

\usepackage{tikz}

\theoremstyle{plain}
\newtheorem{theorem}{Theorem}
\newtheorem{proposition}{Proposition}
\newtheorem*{corollary}{Corollary}
\newtheorem{lemma}{Lemma}

\theoremstyle{definition}
\newtheorem{definition}{Definition}

\newtheorem*{remark}{Remark}

\DeclareMathOperator{\diam}{diam}
\DeclareMathOperator{\loc}{loc}

\DeclareMathOperator{\kernel}{ker}
\DeclareMathOperator{\Image}{Im}

\DeclareMathOperator{\dist}{dist}

\begin{document}

\title{Sharp geometric rigidity of isometries on Heisenberg group}
\author{Daria Isangulova\thanks{The research was carried out within the framework of a state assignment of the Ministry of Education and Science of the Russian Federation for the Institute of Mathematics of the Siberian Branch of the Russian Academy of Sciences (project no. FWNF--2022--0006).}}
\date{}
\maketitle
\begin{abstract}
We prove quantitative stability of isometries on the first
Heisenberg group with sub-Riemannian geometry:  every
$ (1+ \varepsilon)$-quasi-isometry of the John domain of the Heisenberg group
$ \mathbb {H} $ is close to some isometry with order
of closeness  $ \sqrt{\varepsilon} + \varepsilon $ in the uniform norm and with order of closeness $ \varepsilon $ in the Sobolev norm $L_2^1$. Homogeneous dilations show the asymptotic sharpness of the results.

\vspace{0.5cm} \noindent \textbf{2020 Mathematics Subject Classification:} 22E30, 53C17.

\vspace{0.2cm} \noindent \textbf{Keywords:} quasi-isometry,
Heisenberg group, isometry, coercive estimate
\end{abstract}

\section{Introduction}

Geometric rigidity or, stability of isometries, states that for a deformation~$F$ the distance of $DF$ to a suitably chosen proper rotation $Q \in SO(n)$ is dominated by the distance function of $DF$ to $SO(n)$. Both distances can be measured in different norms. Also we can estimate distance from the mapping $F$ to some isometry.
Geometric rigidity plays a central role in models in nonlinear elasticity, it can be considered as 
a suitable nonlinear version of Korn's inequality.

In 1961 F.~John proved geometric rigidity in $\mathbb R^n$, $n\geqslant 2$ \cite{john}:
\textit{for a~locally $(1+\varepsilon)$-bi-Lipschitz mapping $F\colon U\to \mathbb{R}^n$,
where $U$ is an~open set in $\mathbb{R}^n$ and $\varepsilon<1$, there exist a~rotation~$A$ and a vector $a\in \mathbb R^n$
satisfying}
\begin{equation}
\label{eq:Sob_stab_Eucl}
\|DF-A\|_{p,U}
\leqslant C_1 p\varepsilon |U|^{1/p}
\end{equation}
\textit{and}
\begin{equation}
\label{eq:uni_stab_Eucl}
\sup_{x\in U}|F(x)-(a+Ax)|\leqslant C_2\diam(U) \varepsilon.
\end{equation}
F.~John established \eqref{eq:uni_stab_Eucl} for domain $U$ of a~special kind,
now called a~John domain, and \eqref{eq:Sob_stab_Eucl} on cubes.
Later Yu.~G.~Reshetnyak \cite{resh}
established \eqref{eq:Sob_stab_Eucl}  and \eqref{eq:uni_stab_Eucl}
on John domains without constraints on $\varepsilon$ using a~different method.

Friesecke, James and M\"uller showed geometric rigidity in $L^2$-norms for arbitrary Sobolev mapping on a bounded Lipschitz domain $U$ \cite{FJM}: for each $F\in W_2^1(U,\mathbb R^n)$ there is an associated rotation $A\in SO(n)$ such that
$$
\|DF-A\|_{2,U}
\leqslant C_3 \|\dist (DF, SO(n))\|_{2,U}.
$$
Recently Conti, Dolzmann and M\"uller established analogous estimate  under mixed
growth conditions \cite{CDM}.

The geometric rigidity problem can be formulated on any manifold with a~notion
of differential whose tangent space carries an~action
of a~``model'' isometry group.
In this paper, we study the geometric rigidity problem on the Heisenberg
group $\mathbb{H}$ with sub-Riemannian metric. 

Heisenberg group $\mathbb H$ is 3-dimensional nilpotent contact Lie group. Quasiconformal analysis on the Heisenberg group was developed by A.~Kor\'anyi and H.~M.~Reimann \cite{kor-rei-85}. Theory of mappings with bounded distortion was built in the papers of S. K. Vodopyanov \cite{vod-99}, L. Capogna \cite{cap}, N. S. Dairbekov \cite{dair-00}, et al. Introduction to Heisenberg group can be found in the book \cite{cap-book}.

N. Arcozzi and D. Morbidelli \cite{arcozzi-morb} investigated the geometric rigidity
problem for locally $(1+\varepsilon)$-bi-Lipschitz mappings of the Heisenberg group $\mathbb H$
following the idea of John's paper \cite{john}. We should note, however, that
the proximity orders ($\varepsilon^{2^{-11}}$ in the uniform norm and $\varepsilon^{2^{-12}}$ in the Sobolev
norm) obtained in \cite{arcozzi-morb}  are obviously far from being optimal. 

S. K. Vodopyanov  and D. V. Isangulova proved sharp quantitative geometric rigidity of isometries on Heisenberg groups $\mathbb H^n$, $n>1$ \cite{vod-isan-dan,vod-isan-MA}. 
The method of their proof follows the Reshetnyak's scheme \cite{resh} of the proof of stability in the Liouville theorem and essentially uses coercive estimates for a differential operator of the first order with constant coefficients whose kernel coincides with the Lie algebra of the isometry group. Unfortunately, on the first Heisenberg group, such an operator has terms with second-order derivatives \cite{I-SMJ-19}, so it's impossible to apply it to the mappings of the Sobolev class $W_{1,\loc}^1$. In this work
we rely on the P.~P.~Belinsky's idea for proving the stability of conformal mappings in $ \mathbb R^n $ in the uniform norm \cite{bel}.

Elements of Heisenberg group can be identified with the elements of  $\mathbb{R}^{3}$
with the following group law:
$$
(x,y,t)\cdot (x',y',t')=(x+x',y+y',t+t'-2xy'+2x'y).
$$
Vector fields
$
X=\frac\partial{\partial x}
+2y\frac\partial{\partial t}$,
$
Y=\frac\partial{\partial y}
-2x\frac\partial{\partial t}$,
$
T=\frac\partial{\partial t}=-\frac{1}{4} [X,Y]
$
form left-invariant basis of the Lie algebra.
Subbundle $H\mathbb{H}$ of tangent bundle spanned by $X$ and $Y$ 
is called \textit{horizontal}. We introduce such a scalar product on $H\mathbb{H}$ that vector fields $X$ and $Y$ are orthonormal.
{\it Carnot -- Carath\'eodory distance}~$d$ is defined as infimum of lengths of all horizontal curves joining two points  (piecewise-smooth curve is  {\it horizontal}, if its tangent vector belogns to  $H\mathbb H$ almost everywhere).

A bounded open proper subset $U$ of
$\mathbb{H}$ with a distinguished
point $\mathbf{x}_*\in U$ is called a
(metric) \textit{John domain} \cite{john,buck_boman} if it satisfies the
following ``twisted cone'' condition:
there exist constants $\beta\geqslant \alpha  > 0$ such that for all
$\mathbf{x} \in U$, there is a curve
$\gamma: [0, l] \to U$
parameterized by arclength such that
$\gamma(0) = \mathbf{x}$,
$\gamma(l) = \mathbf{x}_*$, $l\leqslant \beta$, and
$\operatorname{dist}(\gamma(s),\partial U) \geqslant \frac{\alpha}{l}s$.
The numbers $\alpha$ and $\beta$ are \textit{inner} and \textit{outer radii} of $U$
respectively.

For a Sobolev mapping $F=(f_1,f_2,f_3)\colon U\to\mathbb H$, $U\subset \mathbb H$,  a linear mapping $$D_h F(\mathbf{x})=
\begin{pmatrix}
Xf_1(\mathbf{x})& Y f_1(\mathbf{x})\\
Xf_2(\mathbf{x}) & Yf_2(\mathbf{x})
\end{pmatrix}\colon H_\mathbf x \mathbb H\to H_{F(\mathbf x)}\mathbb H$$ is called an
{\it approximate horizontal differential}.
(Definition of Sobolev mappings is given in Section~\ref{sec:Sobolev}.)

\begin{definition}[\cite{vod-isan-MA}]\label{def:QI}Let $L \geqslant 1$, $U$ be an open set of $\mathbb{H}$, $F \colon U \to
\mathbb{H}$ be a nonconstant mapping of the Sobolev class
$W^1_{1,\loc} (U,\mathbb{H})$. The mapping $F$ is said to be a $L$-{\it quasiisometry} ($F\in QI_L(U)$), 
if $ L^{-1}|\xi|\leqslant |D_h
F(\mathbf{x}) \xi|\leqslant L |\xi| $ for all vectors $\xi\in H_\mathbf{x}
\mathbb{H}$ for almost all $\mathbf{x}\in U$.
\end{definition}

The main result of the paper is the following

\begin{theorem}\label{th:main}
Let $\varepsilon>0$, $U$ be a John domain with inner radius~$\alpha$ and outer
radius~$\beta$ on the Heisenberg group $\mathbb{H}$. Then, for any
$F\in QI_{1+\varepsilon}(U)$, there exists an isometry $\varphi$ for
which
\begin{equation}\label{eq:uni_global}
\sup_{\mathbf{x}\in U}d(F(\mathbf{x}),\varphi(\mathbf{x}))\leqslant
N_1 \frac{\beta^2}{\alpha} (\sqrt{\varepsilon}+\varepsilon)
\end{equation}
and
\begin{equation}\label{eq:Sob_global}
\int\limits_{U}\exp\Bigl(\Bigl(\frac{\beta}{\alpha}\Bigr)^{5}
\frac{N_2|D_hF(\mathbf{x})-D_h\varphi(\mathbf{x})|}{\varepsilon}\Bigr)\,d\mathbf{x}
\leqslant 16 |U|.
\end{equation}
Here the constants $N_1$ and $N_2$ are independent of $U$ and
$F$.
\end{theorem}

The example of dilations $\delta_{1+\varepsilon}(x,y,t)=((1+\varepsilon)x,(1+\varepsilon)y, (1+\varepsilon)^2t)$ shows
that the proximity orders in Theorem \ref{th:main}
are asymptotically sharp.

Theorem 1 entirely closes the problem of sharp geometric rigidity on Heisenberg groups. Geometric rigidity on $\mathbb H^n$, $n>1$, is formulated absolutely analogously, see Theorem 1 in \cite{vod-isan-MA}. 

First, we prove geometric rigidity in the uniform norm, that is inequality~\eqref{eq:uni_global}. 
The proof is based on the local rigidity of isometries for the following pseudometric $d^H\colon \mathbb{H}\times\mathbb{H}\to [0,\infty)$:
$$
d^H((x,y,t),(x',y',t'))=\sqrt{(x-x')^2+(y-y')^2}.
$$
The function $d^H(\cdot,\cdot)$ satisfies all the axioms of the metric
except
that it can vanish on two different points.

\begin{theorem}\label{th:local_d^H}
Let $\varepsilon>0$, $F\in QI_{1+\varepsilon}(B(\mathbf 0,4))$, $B(\mathbf 0,4)=\{\mathbf x\in\mathbb H \mid d(\mathbf x,\mathbf 0)<4\}$. There exists an isometry $\psi$ for
which
$$
\sup_{\mathbf{x}\in B(\mathbf 0,1)}d^H(F(\mathbf{x}),\psi(\mathbf{x}))\leqslant
N_3  \varepsilon.
$$
Here the constant $N_3$ is independent of $F$.
\end{theorem}

We apply Proposition \ref{prop:local->global} from Section 2 in order to pass from local result of Theorem \ref{th:local_d^H} to global inequality \eqref{eq:uni_global}. In addition, we obtain the following 

\begin{corollary}
Let $\varepsilon>0$, $U$ be a John domain with inner radius~$\alpha$ and outer
radius~$\beta$ on the Heisenberg group $\mathbb{H}$, 
$F\in QI_{1+\varepsilon}(U)$. Then
$$
\sup_{\mathbf{x}\in U}d^H(F(\mathbf{x}),\varphi(\mathbf{x}))\leqslant
N_4 \frac{\beta^2}{\alpha} \varepsilon,
$$
where $\varphi$ is the isometry from Theorem {\rm\ref{th:main}}. 
Constant $N_4$ is independent of $U$ and
$F$.
\end{corollary}

To prove inequality \eqref{eq:Sob_global} from Theorem 1 we establish local geometric rigidity in the Sobolev norm $L_2^1$.

\begin{theorem}\label{th:local_Sob}
Let $\varepsilon>0$, $F\in QI_{1+\varepsilon}(B(\mathbf 0,4))$, $\psi$ is an isometry from
Theorem~{\rm\ref{th:local_d^H}}. Then
$$
\|D_h F(\mathbf{x})-D_h\psi(\mathbf{x})\|_{2,B(\mathbf 0,9/10)}\leqslant N_5 \varepsilon.
$$
Here the constant $N_5$ is independent of $F$, $\|\cdot\|_{2,B(\mathbf 0,9/10)}$ is $L_2$-norm in $B(\mathbf 0,9/10)$.
\end{theorem}

The proof of Theorem \ref{th:local_Sob} is based on the coercive estimate for some differential operator (Lemma \ref{lem:coercive}) and self-improving integrability property of mappings with bounded oscillation (Lemma \ref{lem:BSO}). 
Local result of Theorem \ref{th:local_Sob} means that $D_hF$ is a $BMO$ mapping. 
Proposition \ref{prop:Sob_norm:local_to_global} of Section 5 helps us to deduce global inequality \eqref{eq:Sob_global} from local one. It shows that $BMO$ mapping is exponentially integrable. This finishes the proof of Theorem \ref{th:main}.

The work is built as follows. In Section 2 we give definitions and auxiliary results. Sobolev spaces, isometries and quasi-isometries are defined there with necessary properties.
Proposition \ref{prop:local->global} shows how to pass from local to global results both in Carnot--Carath\'eodory metric and pseudometric $d^H$. Section 3 is a preparatory chapter to the proof of Theorem \ref{th:local_d^H}. We introduce there measurement of closeness between quasi-isometry and the group of isometries. Theorem \ref{th:local_d^H} is proved in Section 4.
Section 5 is devoted to rigidity in the Sobolev norm. 

Results of the paper were announced in \cite{I_2019:DAN}.

\medskip
{\sc Acknowledgments}

The author would like to thank Professor Sergei Vodopyanov for setting the problem and permanent strong belief in success, and Professor Pierre Pansu for useful consultations back in 2005.

\section{Definitions and Auxiliary Results}

\textit{Homogeneous norm}
$\rho(x,y,t)=((x^2+y^2)^2+t^2)^{1/4}$ defines a \textit{Heisenberg metric}~$\rho$:
$\rho(\mathbf{a},\mathbf{b})=\rho(\mathbf{a}^{-1}\cdot \mathbf{b})$ for any
$\mathbf{a},\mathbf{b}\in \mathbb{H}$. (Further, we will denote elements of $\mathbb{H}$ by the bold symbols.)

Metrics~$d$ and~$\rho$ are equivalent \cite{NSW}:
there is a constant $c_0>1$ such that
$\frac{1}{c_0}\rho(\mathbf{a},\mathbf{b})\leqslant d(\mathbf{a},\mathbf{b})
\leqslant c_0 \rho(\mathbf{a},\mathbf{b})$
for all $\mathbf{a},\mathbf{b}\in \mathbb{H}$.

For a ball $B(\mathbf{a},r)=\{\mathbf{u}\in \mathbb{H}\colon d(\mathbf{a},\mathbf{u})<r\}$,
we have $|B(\mathbf{a},r)|=r^{4}|B(\mathbf{0},1)|$, where
$|\cdot|$ is the Lebesgue measure on $\mathbb{R}^{3}$. The latter is the
bi-invariant Haar measure on
$\mathbb{H}$.


Sometimes we will use complex notation: $\mathbf x=(z,t)$, $z=x+iy$.
In the complex notation the following vector fields
\begin{equation*}
Z=\frac{1}{2}(X-i\,Y)
   =\frac{\partial}{\partial z}+i\overline{z}\frac{\partial}{\partial t},
   \quad
\overline{Z}=\frac{1}{2}(X +i\,Y)
   =\frac{\partial}{\partial \overline{z}}
   -iz \frac{\partial}{\partial t}
\end{equation*}
constitute left-invariant basis of the horizontal subbundle
$H\mathbb{H}$.

\subsection{Sobolev spaces on Heisenberg group}\label{sec:Sobolev} 
Let
$\Omega$ be a domain in $\mathbb{H}$.
A \textit{Sobolev space}
$W^1_{q}(\Omega)$ ($L_q^1(\Omega)$), $1\leqslant q\leqslant\infty$,
consists of locally-summable functions
$f \colon \Omega\to\mathbb{R}$ having
generalized derivatives
$Xf$ and $Yf$ and
a finite norm (seminorm)
$$
\|f \|_{W^1_q(\Omega)}=
\|f \|_{q,\Omega}+
\|\nabla_{\mathcal{L}}f \|_{q,\Omega}
\quad
\bigl
(\,\|f \|_{L^1_q(\Omega)}=
\|\nabla_{\mathcal{L}}f\|_{q,\Omega}\bigr),
$$
where
$\nabla_{\mathcal{L}}f=(Xf, Yf)\in H\mathbb H$ is a {\it
subgradient} of the function~$f$.
If $f\in W^1_q(U)$ for each open bounded set
$U$, $\overline U\subset\Omega$, then
$f\in W^1_{q,\loc}(\Omega)$.

\begin{definition}
\label{def:Sob-map}
Let
$\Omega$ be a domain in $\mathbb{H}$.
A mapping $F=(f_1,f_2,f_3)\colon \Omega\to \mathbb{H}$
{\it belongs to the Sobolev class} $W_{q,\loc}^1
(\Omega,\mathbb{H})$,
if\\
1) $f_1,f_2$ and $f_3$ are $ACL$ along a.~e. horizontal coordinate integral lines,\\
2) $f_1$ and $f_2$ belong to the Sobolev
class $W_{q,\loc}^1$
and \\
3) vector fields $XF,YF$
belong to the horizontal subbundle $H\mathbb{H}$
almost everywhere.
\end{definition}

The requirement
$XF(\mathbf{x}),YF(\mathbf{x})\in H_{F(\mathbf{x})}\mathbb{H}$ almost everywhere
yields
so called \textit{contact condition}:
\begin{equation}
\label{eq:contact}
Xf_3=2 f_2 Xf_1-2f_1Xf_2
\quad
\text{and}
\quad
Yf_3=2 f_2\, Yf_1-2f_1\,Yf_2.
\end{equation}

Next lemma shows that pseudometric $d^H$ is sufficient to measure distance between Sobolev mappings in the following sense: $d^H$-distance from $F$ to $\Phi$ controls the Heisenberg distance between them. 

\begin{lemma}\label{lem:d^H(F,G)->rho(F,G)}
Let $F,\Phi\in W_{\infty}^1(B(\mathbf{a},r),\mathbb H)$ and
$$
\sup_{\mathbf{x}\in B(\mathbf{a},r)}d^H(F(\mathbf{x}),\Phi(\mathbf{x}))\leqslant \varepsilon r.
$$
Then there is  $\mathbf b=(0,0,\beta)$, $\beta\in \mathbb{R}$, such that
$$
\sup_{\mathbf{x}\in B(\mathbf{a},r)}
\rho(\mathbf b\cdot F(\mathbf{x}),\Phi(\mathbf{x}))
\leqslant \Bigl(\varepsilon+\sqrt{2\varepsilon\|D_hF+D_h\Phi\|_{\infty,B(\mathbf{a},r)}}\Bigr) r.
$$
\end{lemma}

\begin{proof}
Consider $F=(f_1,f_2,f_3)$, $\Phi=(\varphi_1,\varphi_2,\varphi_3)$, and $\mathbf{x}=(x_1,x_2,x_3)\in B(\mathbf{a},r)$. Denote
$\mathbf b=(0,0,\beta)$, $\beta=-((\Phi(\mathbf{a}))^{-1}\cdot F(\mathbf{a}))_3$ that is $(\Phi(\mathbf{a}))^{-1}\cdot F(\mathbf{a})=(*,*,-\beta)$.

For the mapping $G(\mathbf{x})=(\Phi(\mathbf{x}))^{-1}\cdot \mathbf b\cdot F(\mathbf{x})$ we have
$$
G=(g_1,g_2,g_3)=(f_1-\varphi_1,f_2-\varphi_2,f_3+\beta-\varphi_3+2\varphi_1f_2-2\varphi_2f_1).
$$
Contact condition \eqref{eq:contact}
for $F$ and $\Phi$ yields
$$X_ig_3=2(X_if_1+X_i\varphi_1)(f_2-\varphi_2)-2(X_if_2+X_i\varphi_2)(f_1-\varphi_1),\quad i=1,2.
$$
Thus,
$$
\nabla_\mathcal{L} g_3=2(D_h F+D_h\Phi)^tJ\begin{pmatrix}
f_1-\varphi_1\\
f_2-\varphi_2
\end{pmatrix},
\quad
\text{where}
\quad
J=\begin{pmatrix} 0&1\\
-1&0
\end{pmatrix},$$
and 
$|\nabla_\mathcal{L}g_3(\mathbf{x})|\leqslant 2\|D_h F(\mathbf{x})+D_h \Phi(\mathbf{x})\|d^H(F(\mathbf{x}),\Phi(\mathbf{x}))$.

By construction $g_3(\mathbf{a})=0$ and we obviously have
$$
|g_3(\mathbf{x})|\leqslant
\sup_{\mathbf{y}\in B(\mathbf{a},r)}
\{|\nabla_\mathcal{L}g_3 (\mathbf{y})|\}\, d(\mathbf{a},\mathbf{x})\leqslant 2\|D_hF+D_h\Phi\|_{\infty,B(\mathbf{a},r)}\varepsilon r^2.
$$

Finally, we obtain
\begin{multline*}
\rho(\pi_{\mathbf b}\circ F(\mathbf{x}),\Phi(\mathbf{x}))=\rho(G(\mathbf{x}))\leqslant
d^H(F(\mathbf{x}),\mathbf{x})+\sqrt{|g_3(\mathbf{x})|}
\\
\leqslant \Bigl(\varepsilon + \sqrt{2\varepsilon \|D_hF+D_h\Phi\|_{\infty,B(\mathbf{a},r)}}\Bigr) r.
\end{multline*}
\end{proof}

\subsection{Isometries} 

In this subsection it is convenient to consider a point $\mathbf x = (x,y,t)\in\mathbb{H}$
as $\mathbf x=(z,t)$ where $z=x+iy\in \mathbb{C}$.

Every isometry $\Phi$  (denoted as $\Phi\in \mathrm{Isom}$) is the finite composition of the following mappings \cite[p.~35]{rei} (both in Carnot--Carath\'eodory and Heisenberg metrics):
\begin{itemize}
\item[1)] left translation $\pi_\mathbf{a}(\mathbf{x})=\mathbf{a }\cdot \mathbf{x}$, $\mathbf{a}\in
\mathbb{H}$;
\item[2)] rotation $R_A(z,t) = (e^{iA}z,t)$, $A\in [0,2\pi)$;
\item[3)] reflection $\iota(z, t)
=(\overline{z},-t)$.
\end{itemize}

\begin{lemma}\label{lem:sB_isom}
Let $\Phi\in \mathrm{Isom}$ and\ $d^H(\Phi(\mathbf{x}),\mathbf{x})\leqslant r \varepsilon$
for all $\mathbf{x}\in B(\mathbf{a},r)$.
Then
$$
\sup_{\mathbf{y}\in B(\mathbf{a},sr)}d^H(\Phi(\mathbf{y}),\mathbf{y})
\leqslant (2 s +1)r \varepsilon
\quad
\text{for any }s\geqslant 1
$$
and
$$|D_h \Phi(\mathbf{x})-I|\leqslant 2\varepsilon
\quad \text{for all }\mathbf{x}\in \mathbb{H}.$$
\end{lemma}

\begin{proof}
Put $\psi=\delta_{\frac{1}{r}}\circ \pi_{\mathbf{a}^{-1}}\circ \Phi\circ \pi_\mathbf{a} \circ \delta_r$.
Obviously, a mapping $\psi$ is an isometry, $D_h \psi(\mathbf y)=D_h\Phi(\mathbf w)$ and
$d^H(\psi(\mathbf{y}),\mathbf{y})=\frac{1}{r} d^H(\Phi(\mathbf{w}),\mathbf{w})
\leqslant \varepsilon$
for all $\mathbf{y}\in B(\mathbf 0,1)$ (here $\mathbf{w}=\mathbf{a}\cdot \delta_r \mathbf{y}
\in B(\mathbf{a},r)$).

The isometry $\psi$ can be written as $\psi=(\iota\circ) R_B\circ \pi_{\mathbf{b}}$, where
$\mathbf{b}=(b,\beta)$, $b\in \mathbb{C}$, $\beta\in \mathbb{R}$,
$B\in [0,2\pi)$.
It is easy to see that $|b|=d^H (\psi(\mathbf{0}),\mathbf{0})\leqslant \varepsilon$.

1) Suppose $\psi=\iota\circ R_B\circ \pi_{\mathbf{b}}$.
Then $d^H(\psi(z,t),(z,t))=|e^{-iB}\overline{z}+\overline{b}-z|\leqslant \varepsilon$
and, hence, 
$$
|e^{-iB}\overline{z}-z|\leqslant 2\varepsilon \quad
\text{for any } z\in \mathbb{C},\ |z|\leqslant 1.
$$
Therefore
$$
d^H(\psi(\mathbf{y}),\mathbf{y})=d^H(\psi(\delta_s \mathbf{x}),\delta_s\mathbf{x})
=
|e^{-iB}(\overline{b}+s\overline{z})-sz|
\leqslant
\varepsilon+
2s\varepsilon$$
for all $\mathbf{y}=\delta_s \mathbf{x}\in B(\mathbf{0},s)$, $\mathbf{x}=(z,t)\in B(\mathbf{0},1)$.
Since
$$2\varepsilon\geqslant \sup_{|z|\leqslant 1}
|e^{-iB}\overline{z}-z|=
\sup_{z=e^{iA},A\in[0,2\pi)}
|e^{-i(B+2A)}-1|=2
$$
it follows $\varepsilon\geqslant 1$.
We have
$$
D_h\psi=
\begin{pmatrix}
\cos B & -\sin B\\
-\sin B&-\cos B
\end{pmatrix}
\quad
\text{and}
\quad
|D_h \psi-I|=2\leqslant 2\varepsilon.
$$

2) Suppose now $\psi=R_B\circ \pi_{\mathbf{b}}$.
We have 
$$d^H(\psi(1,0),(1,0))=|e^{iB}(b+1)-1|\leqslant \varepsilon.$$
Therefore, 
$$|D_h \psi-I|\equiv |e^{iB}-1|\leqslant |e^{iB}(b+1)-1|+|b|\leqslant 2\varepsilon.$$
Finally,
$$
d^H(\psi(\mathbf{y}),\mathbf{y})=d^H(\psi(\delta_s \mathbf{x}),\delta_s \mathbf{x})
=
|e^{iB}(b+sz)-sz|
\leqslant \varepsilon+2s\varepsilon$$
for all $\mathbf{y}=\delta_s \mathbf{x}\in B(\mathbf{0},s)$, $\mathbf{x}=(z,t)\in B(\mathbf{0},1)$.

Now we return to the initial isometry $\Phi=\pi_\mathbf{a}\circ \delta_r \circ \psi \circ \delta_{\frac{1}{r}} \circ \pi_{\mathbf{a}^{-1}}$ and the lemma follows.
\end{proof}

\subsection{Quasiisometries. Local qualitative rigidity}

Every quasiisometry $F$ of the class $QI_L(U)$
is locally $L$-Lipschitz on $U$ \cite[Lemma 1]{vod-isan-MA}. If in
addition $F$ is a local homeomorphism then $F$ is locally $L$-bi-Lipschitz.
Conversely, every locally $L$-bi-Lipschitz mapping of an open set $U$ belongs
to $QI_L(U)$.

If~$L$ equals~1 then~$F$ is an isometry.
Obviously, each mapping of the class $QI_L(U)$
is a Sobolev mapping of the class $W_{p,\loc}^1(U, \mathbb{H})$
for any $p\in [1,\infty]$.

Lemma \ref{lem:d^H(F,G)->rho(F,G)} for quasiisometry $F$ and identity mapping $\Phi$ is the following:

\begin{lemma}\label{lem:d^H->rho}
Let $F\in QI_L(B(\mathbf{a},r))$ and
$$
\sup_{\mathbf{x}\in B(\mathbf{a},r)}d^H(F(\mathbf{x}),\mathbf{x})\leqslant \varepsilon r.
$$
Then there is  $\mathbf b=(0,0,\beta)$, $\beta\in \mathbb{R}$, such that
$$
\sup_{\mathbf{x}\in B(\mathbf{a},r)}
\rho(\pi_{\mathbf b}\circ F(\mathbf{x}),\mathbf{x})
\leqslant (\varepsilon+\sqrt{2(L+1)\varepsilon}) r.
$$
\end{lemma}

Next lemma asserts qualitative local rigidity both in the uniform and Sobolev norms. This result is valid for all Heisenberg groups $\mathbb H^n$, $n\geqslant 1$ \cite{vod-isan-MA}.

\begin{lemma}[\mbox{\cite[Lemma 6]{vod-isan-MA}}]
\label{lem:Local_Sob_stab}
For every
$q\in (0,1)$,
there exist nondecreasing functions
$\mu_i(\cdot,q)\colon [0,\infty) \to [0,\infty)$,
$i=1,2$,
such that

{\rm(1)} $\mu_i(t,q)\to 0$ as $t\to 0$,
$i=1,2${\rm;}

{\rm(2)} for each mapping
$f$ of class $QI_{1+\varepsilon}(B(\mathbf{0},1))$, where $B(\mathbf{0},1)\subset \mathbb{H}$,
there exists an~isometry~$\Phi$
satisfying
\begin{align*}
\rho( F(\mathbf x),\Phi(\mathbf x))
&\leqslant q\,\mu_1(\varepsilon ,q)
\quad
\text{for all }\mathbf x\in B(\mathbf 0,q),
\\
\|D_h F-
D_h\Phi\|_{2,B(\mathbf 0,q)}
&\leqslant
|B(\mathbf 0,q)|^{1/2} \mu_2 (\varepsilon,q).
\end{align*}
\end{lemma}

\begin{remark}
For the smaller class of bi-Lipschitz mappings,
Lemma \ref{lem:Local_Sob_stab} is an easy consequence of the results of the paper \cite{arcozzi-morb} where the quantitative estimates for $\mu_i$ is given. In \cite{arcozzi-morb} it is shown that $\mu_1=C\varepsilon^{2^{-11}}$ and $\mu_2=C\varepsilon^{2^{-12}}$.
\end{remark}

\subsection{$KR$-orientation}

Consider Sobolev mapping $F\colon \Omega\to \mathbb H$.
The horizontal differential $D_h F$ is defined almost everywhere and  generates a~morphism $Df$ of graded Lie algebras \cite{vod-07}.

Since
$Df$ is a~homomorphism of graded Lie algebras,
it follows that for almost every $\mathbf x\in \Omega$ there exists
a number $\eta(\mathbf x,F)$ such that
$$
DF(\mathbf x)T=\eta(\mathbf x,F) T.
$$
Furthermore \cite{rei}, $\eta(\mathbf x,F)=\det D_h F(\mathbf x)$
and $\eta(\mathbf x,F)^{2}=\det DF(\mathbf x)$.
Consequently, there are no Sobolev mappings
changing the topological orientation. We now give the definition of
orientation introduced by  A.~Kor\'{a}nyi and H.~M.~Reimann in~\cite{rei}. This orientation considers the sign of the determinant of horizontal differential $D_hF$ instead of the whole differential $DF$.

\begin{definition}
A mapping $F$
of the Sobolev class
$W_{1,\loc}^1(\Omega,\Bbb{H})$
{\it preserves} ({\it changes})
$KR$-{\it orientation} if
$\det D_h F(\mathbf x)>0$ ($\det D_h F(\mathbf x)<0$) for almost all
$\mathbf x\in \Omega$.
\end{definition}

\begin{lemma}[\mbox{\cite[Proposition 3]{isan-smj2}}]\label{lem:KR-orientation}
There is a number $\varepsilon_0 > 0$ such that each quasi-isometry $F\in QI_{1+\varepsilon}(\Omega)$ 
 on a connected open set $\Omega \subset \mathbb H$ with $\varepsilon<\varepsilon_0$ either preserves
$KR$-orientation on the whole $\Omega$ or changes $KR$-orientation on $\Omega$.
\end{lemma}

In \cite{isan-smj2} it is proved for mappings with bounded distortion. We can use this result because quasi-isometries are mappings with bounded distortion. 

\begin{lemma}\label{lem:QI->KR}
Suppose $F\in QI_{1+\varepsilon}(B(\mathbf 0,3/2))$ and 
$$\sup\limits_{\mathbf x\in B(\mathbf 0,1)}d^H(F(\mathbf x),\mathbf x)\leqslant 1/4.$$
There is a number $\varepsilon_1>0$ such that $F$ preserves $KR$-orientation provided $\varepsilon<\varepsilon_1$.
\end{lemma}

\begin{proof}
By Lemma \ref{lem:KR-orientation} quasi-isometry $F$ either preserves $KR$-orientation on the ball $B(\mathbf 0,3/2)$ or reverses $KR$-orientation on the whole ball if $\varepsilon<\varepsilon_0$.
 
In view of Lemma \ref{lem:Local_Sob_stab} there is an isometry $\Phi$ such that
\begin{align*}
\rho( F(\mathbf x),\Phi(\mathbf x))
&\leqslant \mu_1(\varepsilon ,2/3)
\quad
\text{for all }\mathbf x\in B(\mathbf 0,1),
\\
\|D_h F-
D_h\Phi\|_{2,B(\mathbf 0,1)}
&\leqslant
|B(\mathbf 0,1)|^{1/2} \mu_2 (\varepsilon,2/3).
\end{align*}
Take $\varepsilon_1<\varepsilon_0$ such that $\mu_1(\varepsilon_1,2/3)\leqslant 1/4$ and 
$\mu_2(\varepsilon_1,2/3)\leqslant 1$.
We have
$$
d^H(\Phi(\mathbf x),\mathbf x)\leqslant d^H(F(\mathbf x),\Phi(\mathbf x))+d^H(F(\mathbf x),\mathbf x)\leqslant 
\mu_1(\varepsilon ,2/3)+1/4\leqslant 1/2
$$
for all $\mathbf x\in B(\mathbf 0,1)$ if $\varepsilon<\varepsilon_1$.
Lemma \ref{lem:sB_isom} yields 
$$
|D_h \Phi-I|\leqslant 1.
$$
Therefore, $\Phi$ preserves $KR$-orientation. (Here we use an easy fact that $|A-I|>1$ for every $2\times 2$-matrix with $\det A<0$.)

Suppose $F$ reverses $KR$-orientation on $B(\mathbf 0,1)$. Then $\Phi^{-1}\circ F$ also reverses $KR$-orientation and 
$$
|D_h(\Phi^{-1}\circ F)(\mathbf x)-I|>1 \quad\text{for almost all }\mathbf x\in B(\mathbf 0,1).
$$
We get a contradiction since
\begin{multline*}
\|D_h(\Phi^{-1}\circ F)-I\|_{2,B(\mathbf 0,1)}=\|D_hF-D_h\Phi \|_{2,B(\mathbf 0,1)}
\\
\leqslant |B(\mathbf 0,1)|
^{1/2}\mu_2(\varepsilon,2/3)\leqslant |B(\mathbf 0,1)|^{1/2}
\end{multline*}
if $\varepsilon<\varepsilon_1$.
\end{proof}

\subsection{John domains}

In this subsection we demonstrate how to pass from local stability
on balls to stability on John domains for the uniform norm (see Proposition \ref{prop:local->global}).
The analog of this result for functions on Carnot groups can be found in \cite[Lemma 4.2]{Lu}. 

For a domain $U\subset \mathbb{H}$,
denote a distance from a point $\mathbf{x}\in U$ to a boundary $\partial U$
by $d_U(\mathbf{x})=\dist(\mathbf{x},\partial U)=\inf\{d(\mathbf{x},\mathbf{y}):\mathbf{y}\in \partial U\}$.

The following lemma is an easy technical exercise.
\begin{lemma}[\mbox{\cite[Lemma 6]{isan-vod-EMJ}}]\label{lem:chain_balls}
Suppose $\varkappa\geqslant 1$,
$U$ is a John domain in $\mathbb{H}$ with inner radius $\alpha$
and outer radius $\beta$.
Then for each point $\mathbf{x}\in U$
there is a chain of balls
$B_0,\dots,B_k$
satisfying the following
conditions{\rm :}

$(1)$ $B_i=B(\mathbf{x}_i,r_i)$,
$\varkappa B_i=B(\mathbf{x}_i,\varkappa r_i)\subset U$,
$i=0,\dots,k$,
$\mathbf{x}_0=\mathbf{x}_*$ and $\mathbf{x}_k=\mathbf{x}$,
and $d(\mathbf{x},\mathbf{x}_i)\leqslant \frac{\varkappa \beta}{\alpha}r_i$, $i=0,\dots,k-1${\rm ;}

$(2)$ $\frac{2\varkappa-1}{2\varkappa+1}r_{i+1}\leqslant
r_i\leqslant \frac{2\varkappa+1}{2\varkappa-1} r_{i+1}$,
$i=0,\dots,k-1$, and
$\sum\limits_{i=0}^{k-2}r_i+r_k\leqslant 2\beta
${\rm ;}

$(3)$ there is a ball $D_i=B(\mathbf{y}_i,\rho_i)\subset B_i\cap B_{i+1}$ with
$\rho_i=\frac{1}{2}\min\{r_i,r_{i+1}\}$,
$\mathbf{y}_i\in \frac{1}{2}B_i\cap \frac{1}{2}B_{i+1}$ and
$B_k\subset (3+2\frac{(\varkappa+1)\beta}{\alpha})D_i$
for all $i=0,\dots,k-1$.
\end{lemma}

\begin{proposition}\label{prop:local->global}
Suppose $U$ is a John domain in $\mathbb{H}$ with inner radius
$\alpha$, outer radius $\beta$ and distinguished point $\mathbf x_*$, $\varkappa\geqslant 1$, $\sigma>0$. Consider a
mapping $F\in QI_L( U)$ such that for each ball
$B=B(\mathbf{a},r)$, $B(\mathbf{a},\varkappa r)\subset U$, there is
an isometry $\Phi_B$ meeting the condition
$$
\sup_{\mathbf{x}\in B}
d^H(\Phi_B\circ F(\mathbf{x}),\mathbf{x})\leqslant \sigma r.
$$
Consider an isometry $\Phi_*=\pi_{(0,0,s)}\circ \Phi_{B_0}$
with 
$B_0=B(\mathbf x_*,d_U(\mathbf x_*)/\varkappa)$ and $\Phi_*(\mathbf x_*)=(*,*,0)$
$(s=-(\mathbf{x}_*^{-1}\cdot (\Phi_{B_0}\circ F(\mathbf{x}_*)))_3)$.
Then 
\begin{equation}\label{eq:d^H(phi_0(F),id)_Lemma}
\sup_{\mathbf{x}\in U}
d^H(\Phi_*\circ F(\mathbf{x}),\mathbf{x})\leqslant c_1 \sigma \beta
\end{equation}
and
\begin{equation}\label{eq:rho(phi_0(F),id)_Lemma}
\sup_{\mathbf{x}\in U}
\rho(\Phi_*\circ F(\mathbf{x}),\mathbf{x})\leqslant  (c_1\sigma + \sqrt{2(L+1)c_1\sigma}) \beta
\end{equation}
with constant
$c_1=\frac{8\varkappa}{2\varkappa-1} \left(
4L\varkappa\frac{\beta}{\alpha}+2L+1\right)\leqslant 56L\varkappa \frac{\beta}{\alpha}$.
\end{proposition}

\begin{proof}
First we prove \eqref{eq:d^H(phi_0(F),id)_Lemma}.
Fix a point $\mathbf{x}\in U$. Consider a chain of balls
$B_0,\dots,B_k$ from Lemma~\ref{lem:chain_balls}.
Denote $\Phi_i=\Phi_{B_i}$, $G_i=\Phi_i\circ F$, $i=0,\dots,k$,
and $\Phi_{i+1}=\Psi_{i+1}\circ \Phi_i$, $i=0,\dots,k-1$.
Notice that $\Phi_0=\Phi_{B_0}$ is independent of the choice of the point~$\mathbf{x}$.

We have
\begin{multline}\label{eq:d_H(phi_0(F),id)}
d^H(\Phi_0\circ F(\mathbf{x}),\mathbf{x})\leqslant
d^H(\Phi_k\circ F(\mathbf{x}),\mathbf{x})
+
\sum_{i=0}^{k-1}d^H(\Phi_{i+1}\circ F(\mathbf{x}),
\Phi_i\circ F(\mathbf{x}))
\\
\leqslant
\sigma r_k+
\sum_{i=0}^{k-1}d^H(\Psi_{i+1}(G_i(\mathbf{x})),
G_i(\mathbf{x})).
\end{multline}
Estimate $d^H(\Psi_{i+1}(G_i(\mathbf{x})),G_i(\mathbf{x}))$,
$i=0,\dots,k-1$.
Since
\begin{multline*}
d^H(\Psi_{i+1}(\mathbf{y}),\mathbf{y})\leqslant
d^H(\Psi_{i+1}(G_i(\mathbf{y})),\mathbf{y})
+
d^H(\Psi_{i+1}(G_i(\mathbf{y})),\Psi_{i+1}(\mathbf{y}))
\\
\leqslant
(r_i+r_{i+1})\sigma\leqslant
2\sigma \Bigl(1+\frac{2\varkappa+1}{2\varkappa-1}\Bigr)\rho_i
=
\frac{8\varkappa}{2\varkappa-1}\sigma \rho_i
\end{multline*}
for all $\mathbf{y}\in D_i=B(\mathbf{y}_i,\rho_i)$,
Lemma \ref{lem:sB_isom} yields $|D_h \Psi_{i+1}-I|\leqslant \frac{16\varkappa}{2\varkappa-1} \sigma$.
Denote by $g_i\colon U\to \mathbb{C}$ the first complex coordinate
of the mapping $G_i$,
and by $\psi_i\colon \mathbb{C} \to \mathbb{C}$
the first complex coordinate
of the isometry $\Psi_i$. Then
\begin{align*}
d&^H(\Psi_{i+1}(G_i(\mathbf{x})),G_i(\mathbf{x}))
\\
&\leqslant
|\psi_{i+1}(g_i(\mathbf{x}))-\psi_{i+1}(g_i(\mathbf{y}_i))-g_i(\mathbf{x})+g_i(\mathbf{y}_i)|
+
d^H(\Psi_{i+1}(G_i(\mathbf{y}_i)),G_i(\mathbf{y}_i))
\\
&\leqslant
|D_h\Psi_{i+1}-I|d^H(G_i(\mathbf{x}),G_i(\mathbf{y}_i))+
d^H(\Psi_{i+1}(G_i(\mathbf{y}_i)),\mathbf{y}_i)
+
d^H(G_i(\mathbf{y}_i),\mathbf{y}_i)
\\
&\leqslant
\frac{16\varkappa}{2\varkappa-1}\sigma L d(\mathbf{x},\mathbf{y}_i)+\sigma(r_i+r_{i+1})
\quad \text{for }i=0,\dots,k-1.
\end{align*}
We have 
$$
d(\mathbf{x},\mathbf{y}_i)\leqslant d(\mathbf{x},\mathbf{x}_i)+d(\mathbf{x}_i,\mathbf y_i)
\leqslant \frac{\beta}{\alpha}\varkappa r_i+\frac{r_i}{2},
\quad i=0,\dots,k-2,
$$
and
$$
d(\mathbf{x},\mathbf{y}_{k-1})\leqslant \frac{r_k}{2}.
$$
Therefore,
\begin{align*}
d^H(\Psi_{i+1}&(G_i(\mathbf{x})),G_i(\mathbf{x}))
\leqslant
\frac{16\varkappa}{2\varkappa-1}\sigma L\Bigl(\frac{\beta}{\alpha}\varkappa r_i+\frac{1}{2}r_i\Bigr)+\frac{4\varkappa}{2\varkappa-1}\sigma r_i
\\
&= \frac{4\varkappa}{2\varkappa-1} \left(
4L\varkappa\frac{\beta}{\alpha}+2L+1\right)\sigma r_i=\frac{c_1}{2}\sigma r_i \quad \text{for }
i=0,\dots,k-2
\end{align*}
and
\begin{multline*}
d^H(\Psi_{k}(G_{k-1}(\mathbf{x})),G_{k-1}(\mathbf{x}))
\leqslant
\frac{8\varkappa}{2\varkappa-1}\sigma L r_k+
\frac{4\varkappa}{2\varkappa-1} \sigma r_k
\\
=
\frac{4\varkappa}{2\varkappa-1} (2L+1)\sigma r_k
\leqslant \left(\frac{c_1}{2}-1\right)\sigma r_k.
\end{multline*}
In view of \eqref{eq:d_H(phi_0(F),id)} and assertion 2 of Lemma \ref{lem:chain_balls}
we obtain
$$
d^H(\Phi_0\circ F(\mathbf{x}),\mathbf{x})
\leqslant
\sigma r_k+\left(\frac{c_1}{2}-1\right)\sigma r_k+
\frac{c_1}{2}\sigma \sum_{i=0}^{k-2}  r_i
\leqslant c_1\sigma\beta.
$$

Set $\Phi_*=\pi_{(0,0,s)}\circ \Phi_0$ where
$s=-(\mathbf{x}_*^{-1}\cdot (\Phi_0\circ F(\mathbf{x}_*)))_3$.
Obviously,
$$
d^H(\Phi_*\circ F(\mathbf{x}),\mathbf{x})=d^H(\Phi_0\circ F(\mathbf{x}),\mathbf{x})
\leqslant c_1\sigma\beta
$$
and \eqref{eq:d^H(phi_0(F),id)_Lemma} follows.

 Now we prove \eqref{eq:rho(phi_0(F),id)_Lemma}. 
As in the proof of Lemma \ref{lem:d^H->rho} denote
$H=(h_1,h_2,h_3)\colon \mathbf{y}\mapsto\mathbf{y}^{-1}\cdot(\Phi_*\circ F(\mathbf{y}))$.
Then 
$$
|\nabla_\mathcal{L}h_3(\mathbf{y})|\leqslant 2(L+1)d^H(\Phi_*\circ F(\mathbf{y}),\mathbf{y})
\leqslant 2(L+1)c_1 \sigma\beta
$$
for all $\mathbf{y}\in U$
and
$$
|h_3(\mathbf{x})|\leqslant \sup_{\mathbf{y}\in\gamma}
\{|\nabla_\mathcal{L}h_3(\mathbf{y})|\}\,
l\leqslant 2(L+1)c_1\sigma \beta^2
$$
for a curve $\gamma$ joining $\mathbf x$ and $\mathbf x_*$ from definition of John domain.
(Here we used $h_3(\mathbf x_*)=0$.)
Finally,
\begin{multline*}
\rho(\Phi_*\circ F(\mathbf x),\mathbf x)=\rho(H(\mathbf{x}))\leqslant d^H(\Phi_*\circ F(\mathbf{x}),\mathbf{x})+
\sqrt{|h_3(\mathbf{x})|}
\\
\leqslant (c_1\sigma + \sqrt{2(L+1)c_1\sigma })\beta.
\end{multline*}
\end{proof}

\section{Measurement of closeness and normalization}

In this section we build background for proving local quantitative rigidity (Theorem \ref{th:local_d^H}).
First of all we introduce a function $\lambda$ which measures the distance between quasi-isometries and isometries with respect to pseudometric $d^H$. Then we construct a normalized mapping $F_N$ which is close in some sense to the ``optimal'' one.

\begin{definition}
For a mapping $F\colon B(\mathbf{0},1)\to \mathbb{H}$
 define
$$
\lambda_F=\inf_{\Phi\in \mathrm{Isom}}
\sup_{\mathbf{x}\in B(\mathbf{0},1)} d^H(\Phi\circ F(\mathbf{x}),\mathbf{x}).
$$
Given $L,\gamma\geqslant 1$ put
$$
\lambda_L^\gamma=
\sup\{\lambda_F\ \mid\
F\in QI_L(B(\mathbf{0},\gamma))\}.
$$
\end{definition}
In the rest of the paper we will
estimate $\lambda_L^4$
in the terms of $L-1$. Denote $\lambda_L=\lambda_L^4$.

\medskip
\noindent{\bf Properties of $\lambda$:}

1) If the set $F(\overline{B(\mathbf{0},1)})$ is bounded then
$\lambda_F$ is obtained on some $\Phi\in \mathrm{Isom}$;

2) $\lambda_L^\gamma\to 0$ as $L\to 1$ ($\lambda_L^\gamma\leqslant \mu_1(L-1,1/\gamma)$ in view of  Lemma \ref{lem:Local_Sob_stab});

3) $\lambda_L^{\gamma_1}\leqslant \lambda_L^{\gamma_2}$ if $\gamma_1\geqslant \gamma_2$;

4) $\lambda_L^{\gamma_2}\leqslant 56 L\frac{\gamma_1}{\gamma_2}
\lambda_L^{\gamma_1}$ if $\gamma_1\geqslant \gamma_2$ by Proposition \ref{prop:local->global} (since a ball $B(\mathbf{a},r)$ is a John domain with $\alpha=\beta=r$);

5) For any mapping $F\in QI_L(U)$,
there is an isometry $\Phi$
such that
$d^H(\Phi\circ F(\mathbf{x}),\mathbf{x})\leqslant r\lambda_L^\gamma$
for all $\mathbf{x}\in B(\mathbf{a},r)$
if $B(\mathbf{a},\gamma r)\subset U$.

\medskip

\begin{lemma}\label{lem:d^H_on_B(0,3/2)}
Let $F\in QI_L(B(\mathbf{0},4))$
and
$$\sup_{\mathbf{x}\in B(\mathbf{0},1)}d^H(F(\mathbf{x}),\mathbf{x})\leqslant c\lambda_L,
\quad c\geqslant 1.
$$
Then
$$\sup_{\mathbf{x}\in B(\mathbf{0},3/2)}d^H(F(\mathbf{x}),\mathbf{x})
\leqslant (7 c+4)  \lambda_L.$$
\end{lemma}

\begin{proof}
Fix a point $\mathbf{x}\in B(\mathbf{0},3/2)\setminus B(\mathbf{0},1)$.
Consider a ball $B(\mathbf{x}_0,1/2)$ containing $
\mathbf x$
such that $\mathbf x_0\in\partial B(\mathbf 0,1)$ ($d(\mathbf{x}_0,\mathbf{0})=1$).
Notice that
$B(\mathbf{x}_0,2)\subset B(\mathbf{0},4)$.
There is an isometry $\Phi$ such that
$$
\sup_{\mathbf{y}\in B(\mathbf{x}_0,1/2)}
d^H(\Phi\circ F(\mathbf{y}),\mathbf{y})\leqslant \frac{1}{2} \lambda_L.
$$

Consider a ball $G$ of radius $1/4$ in the intersection $B(\mathbf{0},1)\cap
B(\mathbf{x}_0,1/2)$.
For $\mathbf{y}\in G$, we have
\begin{multline*}
d^H(\Phi^{-1}(\mathbf{y}),\mathbf{y})\leqslant
d^H(\Phi^{-1}(\mathbf{y}),F(\mathbf{y}))+d^H(F(\mathbf{y}),\mathbf{y})
\\
=
d^H(\Phi\circ F(\mathbf{y}),\mathbf{y})+d^H(F(\mathbf{y}),\mathbf{y})
\leqslant \lambda_L\Bigl(\frac{1}{2}+c\Bigr)
\end{multline*}
In view of Lemma \ref{lem:sB_isom} it follows
$d^H(\Phi^{-1}(\mathbf{y}),\mathbf{y})\leqslant 7 \lambda_L
(\frac{1}{2}+c)$ for all $\mathbf{y}\in 3G$.
We obtain
\begin{multline*}
d^H(F(\mathbf{x}),\mathbf{x})\leqslant
d^H (\Phi\circ F(\mathbf{x}),\mathbf{x})
+d^H(\Phi^{-1}(\mathbf{x}),\mathbf{x})
\\
\leqslant \frac{1}{2}\lambda_L+7\lambda_L \Bigl(\frac{1}{2}+c\Bigr)
\leqslant (7c+4) \lambda_L
\end{multline*}
since $\mathbf{x}\in B(\mathbf{x}_0,1/2)\subset 3G$.
\end{proof}

\begin{lemma}\label{lem:auxiliary}
Suppose $F\in QI_L(B(\mathbf{0},4))$ and
$\sup\limits_{\mathbf{x}\in B(\mathbf{0},3/2)}d^H(F(\mathbf{x}),\mathbf{x})\leqslant c\lambda_L$.
Consider a number $r$, $0<r<3/8$, and a point  $\mathbf x_0$ satisfying $B(\mathbf{x}_0,4r)\subset B(\mathbf{0},3/2)$.
Then
$$d^H(\mathbf{x}^{-1}\cdot F(\mathbf{x}),\mathbf{x}_0^{-1}\cdot F(\mathbf{x}_0))
\leqslant 2\Bigl(1+\frac{(r+c)}{r^2}d^H(F(\mathbf x),F(\mathbf x_0))\Bigr) r \lambda_L$$
for every $\mathbf x\in B(\mathbf x_0,r)$.
\end{lemma}

\begin{proof}
Consider a ball $B_0=B(\mathbf{x}_0,r)$, $B(\mathbf{x}_0,4r)\subset B(\mathbf{0},3/2)$.
Then there is an isometry $\Phi$
such that 
$$\sup_{\mathbf{y}\in B_0}
d^H(\Phi\circ F(\mathbf{y}),\mathbf{y})=r \lambda_F\leqslant r\lambda_L.$$
We have
$$
d^H(\Phi(\mathbf{y}),\mathbf{y})\leqslant
d^H(\Phi\circ F(\mathbf{y}),\mathbf{y})+d^H(\Phi\circ F(\mathbf{y}),\Phi(\mathbf{y}))
\leqslant (r+c)\lambda_L,\quad
\mathbf{y}\in B_0.
$$
In view of Lemma \ref{lem:sB_isom} it follows
$|D_h \Phi -I|\leqslant 2 \frac{r+c}{r}\lambda_L$.

Therefore,
\begin{multline*}
d^H(\mathbf{x}^{-1}F(\mathbf{x}),\mathbf{x}_0^{-1}F(\mathbf{x}_0))
\leqslant
d^H(\Phi\circ F(\mathbf{x}),\mathbf{x})
+d^H(\Phi\circ F(\mathbf{x}_0),\mathbf{x}_0)
\\
+d^H(F(\mathbf{x})^{-1}\cdot \Phi(F(\mathbf{x})),F(\mathbf{x}_0)^{-1}
\cdot\Phi(F(\mathbf{x}_0)))
\\
\leqslant
2 r \lambda_L+|D_h \Phi-I|\,d^H(F(\mathbf{x}),F(\mathbf{x}_0))
\\
\leqslant
2 r\lambda_L\Bigl(1+\frac{r+c}{r^2}d^H(F(\mathbf{x}),F(\mathbf{x}_0))\Bigr)
\end{multline*}
for every $\mathbf x\in B(\mathbf x_0,r)$.
\end{proof}

Consider now a point $\mathbf x\in\mathbb H$ in complex notation: $\mathbf x=(z,t)$, $z\in\mathbb C$,
$t\in\mathbb R$.

\begin{definition}
Suppose $F\colon \overline{B(\mathbf{0},1)}\to \mathbb{H}$ is a
continuous mapping,
$f\colon \overline{B(\mathbf{0},1)}\to \mathbb{C}$ is the first complex
coordinate function of $F$.
Set
$\mathbf{a}=F(\mathbf{0})^{-1}$ and $A=-\arg (f(1,0)-f(0,0))$.

The mapping $F_N=R_A\circ \pi_\mathbf{a}\circ F$
is a \textit{normalization} of the mapping $F$.
\end{definition}

\medskip
\noindent{\bf Properties of normalization:}

1) $F_N(\mathbf 0)=\mathbf 0$;

2) $F_N(1,0)=(b,\beta)$, $\mathrm{Re}\, b\geqslant 0$, $\mathrm{Im}\, b=0$,  $\beta\in\mathbb R$;

3) If $\Psi\circ F(\mathbf 0)=\mathbf 0$, $\Psi\circ F(1,0)=(b,\beta)$, $\mathrm{Re}\, b\geqslant 0$, $\mathrm{Im}\, b=0$,  $\beta\in\mathbb R$, for  some isometry $\Psi$, then either $\Psi\circ F=F_N$ or $\Psi\circ F=\iota\circ F_N$ where $\iota$ is reflection, $\iota(z,t)=(\overline z,-t)$.

\medskip

\begin{lemma}\label{lem:d^H(F_N,id)}
There is a number $\varepsilon_2>0$ such that
$$
\sup\limits_{\mathbf{x}\in B(\mathbf{0},1)}d^H(F_N(\mathbf{x}),\mathbf{x})\leqslant 6 \lambda_F
$$
for every $F\in QI_{1+\varepsilon}(B(\mathbf 0,3/2))$ provided $\varepsilon<\varepsilon_2$ and
$F$ preserves $KR$-orientation.
\end{lemma}

\begin{proof} 

Assume $\lambda_F=\sup\limits_{\mathbf{x}\in B(\mathbf{0},1)}d^H(\Psi\circ F(\mathbf{x}),\mathbf{x})$, $\Psi\in\operatorname{Isom}$,
and put $G=\Psi\circ F$.
Let $G_N=\Phi\circ G$. 
We have $\Phi=R_B\circ \pi_\mathbf{b}$
with $\mathbf{b}=(b,\beta)=G(\mathbf{0})^{-1}$ and
$g(1,0)-g(0,0)=r e^{-iB}$. (Here  $g$ is the first complex coordinate function of the mapping $G$.)

Obviously, $|b|\leqslant \lambda_F$
and $|r e^{-iB}-1|\leqslant 2\lambda_F$.
It follows
$$|e^{-iB}-1|\leqslant |r e^{-iB}-1|+|r-1|
\leqslant 2 |
r e^{-iB}-1|\leqslant 4\lambda_F.
$$
From here $d^H(\Phi(\mathbf{x}),\mathbf{x})=|e^{iB}(b+z)-z|
\leqslant |b|+|e^{iB}-1|\leqslant 5\lambda_F$
for any $\mathbf{x}=(z,t)\in B(\mathbf{0},1)$.
Thus, for any $\mathbf{x}\in B(\mathbf{0},1)$ we have
$$
d^H(G_N(\mathbf{x}),\mathbf{x})\leqslant d^H(\Phi\circ G(\mathbf{x}),\Phi(\mathbf{x}))
+
d^H(\Phi(\mathbf{x}),\mathbf{x})\leqslant \lambda_F+5\lambda_F=6\lambda_F.
$$
In view of the property 3 of normalization, either $G_N=F_N$ or $G_N=\iota\circ F_N$. 

Suppose our lemma is incorrect and
$$\sup\limits_{\mathbf{x}\in B(\mathbf{0},1)}d^H(\iota\circ F_N(\mathbf{x}),\mathbf{x})\leqslant 6 \lambda_F.$$

Recall
$\lambda_F\leqslant \lambda_{1+\varepsilon}^{3/2}\to 0$ as $\varepsilon\to 0$. There is $\varepsilon_2<\varepsilon_1$ such that $\lambda_F<1/24$ for $F\in QI_{1+\varepsilon}(B(\mathbf 0,3/2))$, $\varepsilon<\varepsilon_2$. Therefore, by Lemma \ref{lem:QI->KR} $\iota\circ F_N$ preserves $KR$-orientation. We come to a contradiction, since $\det D_hF=\det D_h F_N>0$ and $\det D_h\iota<0$. 
 \end{proof}

\section{Proof of Theorem \ref{th:local_d^H}}

To prove Theorem \ref{th:local_d^H} it suffices to show that $\lambda_L\leqslant C(L-1)$ as $L$ tends to 1.
Recall
$$
\lambda_L=
\sup\{\lambda_F\ \mid\
F\in QI_L(B(\mathbf{0},4))\},
$$
where
$$
\lambda_F=\inf_{\Phi\in \mathrm{Isom}}
\sup_{\mathbf{x}\in B(\mathbf{0},1)} d^H(\Phi\circ F(\mathbf{x}),\mathbf{x}).
$$

Consider $G\in QI_L(B(\mathbf{0},4))$ with $\lambda_G=\sup\limits_{\mathbf{x}\in B(\mathbf{0},1)} d^H(G(\mathbf{x}),\mathbf{x})\geqslant 9\lambda_L/10$. 


\subsection{Step 1}

Let $L^{m}\leqslant \Lambda<L^{m+1}$, $m\in\mathbb N$. A number $\Lambda<\min\{1+\varepsilon_2,3/2\}$ is such that
$\lambda_\Lambda\leqslant M$. Quantity $M<1/4$ we will choose later. 

We have $$\sup\limits_{\mathbf{x}\in B(\mathbf{0},1)} d^H(G(\mathbf{x}),\mathbf{x})=\lambda_G\leqslant \lambda_L\leqslant \lambda_\Lambda<1/4\quad \text{and}\quad L<1+\varepsilon_2<1+\varepsilon_1.$$ Then
by Lemma \ref{lem:QI->KR} $G$ preserves $KR$-orientation.  Denote $F=G_N$. 
Lemma \ref{lem:d^H(F_N,id)} implies
\begin{equation}\label{eq:d^H(F(x),x)}
\frac{9}{10}\lambda_L\leqslant \sup_{\mathbf{x}\in B(\mathbf{0},1)}
d^H(F(\mathbf{x}),\mathbf{x})\leqslant 6\lambda_L.
\end{equation}

Lemma \ref{lem:d^H_on_B(0,3/2)}
implies
\begin{equation}\label{eq:D^H(F(x),x)_on_B(0,3/2)}
\sup_{\mathbf{x}\in B(\mathbf{0},3/2)}d^H(F(\mathbf{x}),\mathbf{x})\leqslant  c_3 \lambda_L
\quad \text{with }c_3=46.
\end{equation}

We have $F^l$ is locally $L^l$-Lipschitz and $F^l(\mathbf 0)=\mathbf 0$. Therefore, 
$$
d(F^l(\mathbf x),\mathbf 0)=d(F^l(\mathbf x),F^l(\mathbf 0))\leqslant L^l d(\mathbf x,\mathbf 0)<\Lambda d(\mathbf x,\mathbf 0)<\frac{3}{2}d(\mathbf x,\mathbf 0)
$$ 
for all $l\leqslant m$ and every $\mathbf x\in \overline{B(\mathbf{0},1)}$.
It follows
$F^l(\overline{B(\mathbf{0},1)})\subset B(\mathbf{0},3/2)$. 

Take $\mathbf{x}\in \overline{B(\mathbf{0},1)}$ and $l\leqslant m$, $l\in\mathbb N$. 
Inequality \eqref{eq:D^H(F(x),x)_on_B(0,3/2)} implies
\begin{equation}\label{eq:d^H(F^l(x),x)}
d^H(F^{l}(\mathbf{x}),\mathbf{x})\leqslant
d^H(F^l(\mathbf{x}),F^{l-1}(\mathbf{x}))+\dots+
d^H(F(\mathbf{x}),\mathbf{x})\leqslant c_3 l\lambda_L
\end{equation}
and
\begin{equation*}
d^H(F^{l+1}(\mathbf{x}),F(\mathbf{x}))\leqslant
d^H(F^{l+1}(\mathbf{x}),F^{l}(\mathbf{x}))+\dots+
d^H(F^2(\mathbf{x}),F(\mathbf{x}))\leqslant c_3 l\lambda_L.
\end{equation*}

Lemma \ref{lem:d^H->rho} yields
$$
d(F^l(\mathbf{x}),\mathbf{x})\leqslant c_0(c_3 l\lambda_L+\sqrt{2(L^l+1)c_3 l\lambda_L})
$$
where $c_0$ is the coefficient of equivalency between Heisenberg metric $\rho$ and Carnot--Carath\'eodory metric $d$.
If $c_3l\lambda_L<1/4$ and $L^l<\Lambda<1.1$ then
$$(c_3 l\lambda_L+\sqrt{2(L^l+1)c_3 l\lambda_L}) \leqslant 2\sqrt{c_3 l\lambda_L}.
$$

Fix some $r$, $0<r<3/8$. Suppose
\begin{equation}\label{eq:D<r^2}
2c_0\sqrt{c_3 l\lambda_L}\leqslant r.
\end{equation}
Then
$$
d(F^l(\mathbf{x}),\mathbf{x})\leqslant r
\quad
\text{and}
\quad
d^H(F^{l+1}(\mathbf{x}),F(\mathbf{x}))\leqslant \frac{r^2}{4c_0^2}.
$$

By Lemma \ref{lem:auxiliary} we have 
\begin{multline}\label{eq:d^H(f^l+1-f^l+f-x)}
d^H\bigl((F^l(\mathbf{x}))^{-1}\cdot F^{l+1}(\mathbf{x}),\mathbf{x}^{-1}\cdot F(\mathbf{x})\bigr)
\\
\leqslant 2\Bigl(1+\frac{(r+c_3)}{r^2}d^H(F^{l+1}(\mathbf{x}),F(\mathbf{x}))\Bigr)r\lambda_L
\leqslant N r \lambda_L
\end{multline}
where $N= 2(1+(1+c_3)/4)=25.5$.

\subsection{Step 2}

Fix $l<m$, $l\in\mathbb N$. Denote by $f_l\colon B(\mathbf 0,3/2)\to\mathbb C$ the first complex coordinate function of
$F^l$, $l=1,\dots,m$. Set $f_{l}(1,0)=r_l e^{i\theta_l}$.
In this subsection we will estimate $|e^{i\theta_{l+1}}-e^{i\theta_l}|$. 

Recall $c_3 l\lambda_L\leqslant \frac{r^2}{4c_0^2}<\frac{1}{4}$ by \eqref{eq:D<r^2}.
Relation \eqref{eq:d^H(F^l(x),x)} at the point $(1,0)$
gives 
$$
|r_l e^{i\theta_l}-1|\leqslant c_3 l \lambda_L<1/4.
$$
This yields
$$
|r_l-1|\leqslant |r_l e^{i\theta_l}-1|\leqslant c_3 l \lambda_L<1/4
$$
and
$$
|e^{i\theta_l}-1|\leqslant |r_l e^{i\theta_l}-1|+|r_l-1|\leqslant 2c_3 l \lambda_L<1/2.
$$
Since $F_N=F$ we have $\theta_1=0$ and $f_1(1,0)=r_1$. In view of \eqref{eq:d^H(F(x),x)} it follows $|r_1-1|\leqslant 6 \lambda_L$.

Write down equation \eqref{eq:d^H(f^l+1-f^l+f-x)}
at the point $(1,0)$:
\begin{equation}\label{eq:d^H(f^l+1-f^l+f-x),x=(1,0)}
|f_{l+1}(1,0)-f_{l}(1,0)-f(1,0)+1|=
|r_{l+1}e^{i\theta_{l+1}}-r_{l}e^{i\theta_{l}}-r_1+1|\leqslant Nr\lambda_L.
\end{equation}
Obviously,
$$
|\Image (r_{l+1}e^{i\theta_{l+1}}-r_le^{i\theta_l})|
=|r_{l+1}\sin \theta_{l+1}-r_l\sin \theta_l |\leqslant Nr\lambda_L
$$
and
$$
|r_{l+1}\sin \theta_{l+1}|\leqslant Nr\lambda_L+|r_{l}\sin \theta_{l}|
\leqslant \dots\leqslant lNr\lambda_L.
$$

Define 3 points in the complex plane (see Figure 1):
\begin{align*}
A&=f_l(1,0)=r_le^{i\theta_l},\\
B&=f_l(1,0)+f_1(1,0)-1=\overline{r}_l e^{i (\theta_l+\xi_l)},\\
C&=\widehat{r}_le^{i\theta_l} \ \text{such that }BC \bot OA.
\end{align*}

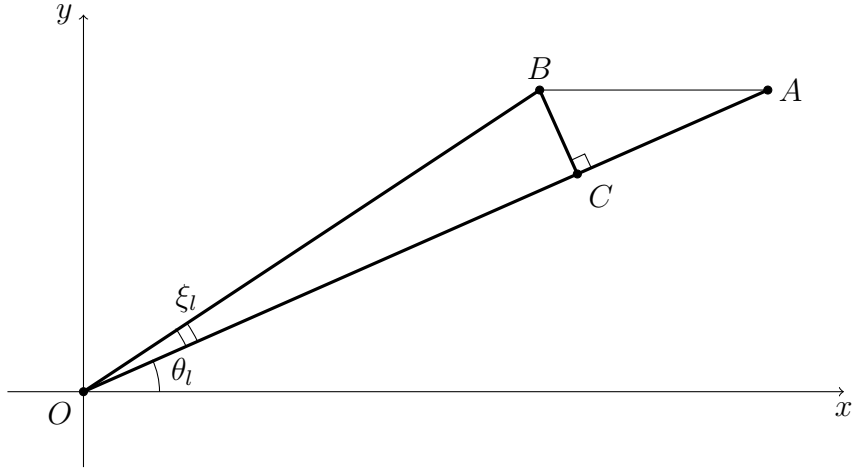
\begin{center}
\begin{tikzpicture}[scale=1]
\draw[->] (-1,0) -- (10,0) node[anchor=north] {$x$};
\draw[->] (0,-1) -- (0,5) node[anchor=east] {$y$};  

\draw [fill,black] (0,0) circle (1.5pt) node[anchor=north east] {$O$};

\draw[very thick] (0,0) -- (9,4);
\draw[fill,black] (9,4) circle (1.5pt) node[anchor=west] {$A$};
\draw (10mm,0) arc (0:25:10mm) node[below=1.5mm,right=1mm] {$\theta_l$};
\draw[thin] (9,4) -- (6,4);

\draw[very thick] (0,0) -- (6,4);
\draw (1.35,0.6) arc (25:35:1.57cm);
\draw (1.5,0.666) arc (25:35:1.64cm)  node[anchor=south] {$\xi_l$};
\draw[fill,black] (6,4) circle (1.5pt) node[anchor=south] {$B$}
(intersection of 0,0 -- 9,4 and 6,4 -- 10,-5) coordinate (t);

\draw[very thick] (6,4) -- (t);
 \draw[thin]  [shift=(t)] [rotate=25] (0mm,2mm) -- (2mm,2mm) -- (2mm,0mm); 
 \draw[fill,black] (t) circle (1.5pt) node[anchor=north west] {$C$};
\end{tikzpicture}
\end{center}

\vspace{-1cm}
\begin{center}{Figure 1.}
\end{center}

We have
$$
|AB|=|f_1(1,0)-1|=|r_1-1|<6\lambda_L,
$$
$$
|OB|=\overline{r}_l\geqslant r_l-|r_1-1|\geqslant 3/4-6\lambda_L
\geqslant 3/4-6M
>1/2 \quad \text{provided }M<1/24.
$$

Our next aim is to estimate $\xi_l$.
Considering triangles $ABC$ and $OBC$ we obtain
$$
|BC|=|AB|\sin \angle CAB=|r_1-1||\sin \theta_l|\leqslant
6\lambda_L \frac{Nr(l-1)\lambda_L}{r_l}
\leqslant
\frac{2}{c_3} N r\lambda_L
$$
and
$$
\sin  \angle COB=|\sin\xi_l|=\frac{|BC|}{|OB|}\leqslant
\frac{2}{c_3 \overline{r}_l} Nr\lambda_L<\frac{4}{c_3} Nr\lambda_L<4M<\frac{1}{2} 
$$
if $M<1/8$.
Therefore, 
\begin{multline}\label{eq:16}
|e^{i\xi_l}-1|=\sqrt{\sin^2 \xi_l+(\cos \xi_l-1)^2}=
\sqrt{2(1-\cos \xi_l)}=2\Bigl|\sin\frac{\xi_l}{2}\Bigr|\\
\leqslant
2|\sin\xi_l|<\frac{8}{c_3} Nr\lambda_L.
\end{multline}

Equation \eqref{eq:d^H(f^l+1-f^l+f-x),x=(1,0)} yields
$$
|r_{l+1}e^{i\theta_{l+1}}-\overline{r}_l e^{i(\theta_l+\xi_l)}|
\leqslant Nr\lambda_L.
$$
It follows
$$
|r_{l+1}-\overline{r}_l|\leqslant
|r_{l+1}e^{i\theta_{l+1}}-\overline{r}_le^{i(\theta_l+\xi_l)}|\leqslant N r\lambda_L
$$
and
\begin{equation}\label{eq:17}
|e^{i \theta_{l+1}}-e^{i(\theta_l+\xi_l)}|\leqslant
\frac{|r_{l+1} e^{i(\theta_{l+1}-\theta_l-\xi_l)}-\overline{r}_{l}|+|r_{l+1}-\overline{r}_l|}{r_{l+1}}
\leqslant
\frac{4}{3} 2 N r\lambda_L.
\end{equation}
Finally, from \eqref{eq:16} and \eqref{eq:17} we derive
\begin{equation}\label{eq:e^itheta_l+1-e^itheta_l}
|e^{i \theta_{l+1}}-e^{i\theta_l}|\leqslant
|e^{i \theta_{l+1}}-e^{i(\theta_l+\xi_l)}|
+|e^{i \xi_l}-1|\leqslant \frac{8}{3} N r\lambda_L
+\frac{8}{c_3} N r\lambda_L<3 N r\lambda_L.
\end{equation}

\subsection{Step 3}

Denote by $g_l\colon B(\mathbf 0,3/2)\to\mathbb C$ the first complex coordinate function of
$(F^l)_N$, $l=1,\dots,m$. 
Since $F^l(\mathbf{0})=\mathbf{0}$ and $f_{l}(1,0)=r_l e^{i\theta_l}$, it yields
$g_l=e^{-i\theta_l}f_l$. Obviously, $g_1=f_1$.

At the point $\mathbf{x}=(z,t)\in \overline{B(\mathbf{0},1)}$ we have
\begin{multline*}
|g_{l+1}(\mathbf{x})-g_l(\mathbf{x})-f_1(\mathbf{x})+z|\leqslant
|f_{l+1}(\mathbf{x})-f_l(\mathbf{x})-f_1(\mathbf{x})+z|
\\
+
|e^{-i\theta_{l+1}} f_{l+1}(\mathbf{x})-e^{-i\theta_l}f_l(\mathbf{x})-f_{l+1}(\mathbf{x})+f_l(\mathbf{x})|
\\
\leqslant
N r\lambda_L+|e^{-i\theta_{l}}-1|\cdot |f_{l+1}(\mathbf{x})-f_l(\mathbf{x})|
+|e^{-i\theta_{l+1}}-e^{-i\theta_l}|\cdot |f_{l+1}(\mathbf{x})|
\\
\leqslant
N r \lambda_L+2c_3 l \lambda_L c_3 \lambda_L+3N r \lambda_L \frac{3}{2}
\\
\leqslant
\Bigl(N+\frac{2c_3 r}{4c_0^2}+3N\frac{3}{2}\Bigr)r\lambda_L
\leqslant
Kr\lambda_L,\quad K=6N=153
\end{multline*}
(we used here \eqref{eq:D^H(F(x),x)_on_B(0,3/2)}, \eqref{eq:D<r^2}, \eqref{eq:d^H(f^l+1-f^l+f-x)},  \eqref{eq:e^itheta_l+1-e^itheta_l} and $r<3/8$).

Summing over integers from 1 to $l$
we obtain
\begin{equation}\label{eq:19}
|g_{l+1}(\mathbf{x})-z-l(f_1(\mathbf{x})-z)|\leqslant Klr\lambda_L=\frac{2}{5}l\lambda_L
\end{equation}
provided we take $r= \frac{2}{5K}$.

Let a point $\mathbf{x}_0=(z_0,t_0)\in \overline{B(\mathbf{0},1)}$
be the point of maximum deviation:
$$|f_1(\mathbf{x}_0)-z_0|=d^H(F(\mathbf{x}_0),\mathbf{x}_0)=\sup\limits_{\mathbf{x}\in \overline{B(\mathbf{0},1)}}d^H(F(\mathbf{x}),\mathbf{x})\geqslant \frac{9}{10}\lambda_L.$$

Take $l\in\mathbb N$, $l<m$. Since $F$ preserves $KR$-orientation, the mapping $F^{l+1}$ also preserves it. $F^{l+1}$ is $L^{l+1}$-quasi-isometry on $B(\mathbf 0,3/2)$ and $L^{l+1}<\Lambda<1+\varepsilon_2$. Therefore, by
Lemma \ref{lem:d^H(F_N,id)}
 and \eqref{eq:19} we deduce 
\begin{multline*}
6\lambda_{L^{l+1}}\geqslant d^H((F^{l+1})_N(\mathbf{x}_0),\mathbf{x}_0)=
|g_{l+1}(\mathbf{x}_0)-z_0|
\\
\geqslant
l|f_1(\mathbf{x}_0)-z_0|-|g_{l+1}(\mathbf{x}_0)-z_0-l(f_1(\mathbf{x}_0)-z_0)|\\
\geqslant
\frac{9}{10}l\lambda_L -\frac{2}{5}l\lambda_L
= \frac{1}{2} l \lambda_L.
\end{multline*}
We immediately obtain
\begin{equation}\label{eq:lambda_L}
l\lambda_L\leqslant 12\lambda_{L^{l+1}}, \quad l=1,\dots,m-1.
\end{equation}

Suppose $m\geqslant 3$. Since $L^m\leqslant \Lambda<L^{m+1}$ it follows
$$
m\geqslant \frac{\ln \Lambda}{\ln L}-1
\quad
\text{and}
\quad
\ln \Lambda -\ln L\geqslant \frac{\ln \Lambda}{2} .
$$
Therefore, \eqref{eq:lambda_L} implies
$$
\lambda_L\leqslant \frac{12\lambda_\Lambda}{\frac{\ln \Lambda}{\ln L}-1}
\leqslant \frac{24\lambda_\Lambda}{\ln \Lambda}\ln L\leqslant C(L-1)
$$
and Theorem 2 is proved provided inequality \eqref{eq:D<r^2} is fulfilled.

\subsection{Step 4}
It rests to prove \eqref{eq:D<r^2}.
Demonstrate by induction that inequality
\eqref{eq:D<r^2} is fulfilled provided $M=\lambda_\Lambda\leqslant \frac{r^2}{96 c_3 c_0^2}$.

Since $\lambda_L$ tends to 0 as $L\to 1$ there is a number $\varepsilon_3>0$
such that $\varepsilon_3<\varepsilon_2$, $\varepsilon_3<0.1$, and $M=\lambda_\Lambda<\frac{r^2}{96 c_3 c_0^2}$ if $\Lambda=1+\varepsilon_3$. 

Base of induction: $c_3\lambda_L\leqslant c_3 M<\frac{\rho^2}{96c_0^2}$.
Let for some $l$, $1\leqslant l<m$, we have $c_3l\lambda_L\leqslant \frac{r^2}{4c_0^2}$.
Then \eqref{eq:lambda_L} implies
$$
l\lambda_L\leqslant 12\lambda_\Lambda\leqslant 12 M\leqslant \frac{12 r^2}{96 c_3 c_0^2}
\quad
\text{and, hence, } c_3 l\lambda_L \leqslant \frac{r^2}{8c_0^2}.
$$
From here
$c_3(l+1)\lambda_L\leqslant \frac{r^2}{4c_0^2}$
and inequality \eqref{eq:D<r^2} holds for $l+1$.
We can proceed till $l=m-1$.

Therefore, Theorem 2 is proved for small $\varepsilon$ (it is sufficient to consider
$\varepsilon<\varepsilon_3/3$). 
Consider the case $L\geqslant 1+\varepsilon_3/3$.
We obviously have
$$
d(F(\mathbf{x}),\mathbf{x})\leqslant
d(F(\mathbf{x}),F(\mathbf{0}))+d(\mathbf{x},\mathbf{0})
\leqslant (L+1)\leqslant \Bigl(1+\frac{6}{\varepsilon_3}\Bigr)(L-1)
$$
for all $\mathbf{x}\in B(\mathbf 0,1)$. 

Theorem 2 is proved.

\section{Rigidity in the Sobolev norm}

The aim of this section is to prove inequality \eqref{eq:Sob_global} from Theorem 1. The main difficulty is to prove local quantitative rigidity in Sobolev norm (Theorem \ref{th:local_Sob}). The proof is based on coercive estimate for special differential operator (Lemma \ref{lem:coercive}) and the connection between this operator and quasi-isometries, see inequality \eqref{eq:main}. The already established quantitative rigidity in the uniform norm also plays an essential role in the proof. To pass from local result of Theorem  \ref{th:local_Sob} to global one in Theorem 1 we use Proposition \ref{prop:Sob_norm:local_to_global}.

\subsection{Coercive estimate}

Let $U$ be a domain in $\mathbb{H}$. Denote by $Q$  the
homogeneous differential operator acting on a mapping $u \colon U
\to \mathbb{C}$ by the following rule:
$$
Qu=\begin{pmatrix}
\frac{1}{2}(Z u + \overline{Zu})\\
\overline{Z}u
\end{pmatrix}.
$$

\begin{lemma}\label{lem:coercive}
Let $p>1$, $B=B(\mathbf 0,1)\subset\mathbb H$. Then there is a constant $C=C(p)>0$ such that
$$
\|D_hu\|_{p,B}\leqslant C \Bigl(\sup\limits_{\mathbf{x}\in
B(\mathbf{0},1)}|u(\mathbf{x})|+\|Qu\|_{p,B}\Bigr)
$$
for every $W_{p}^1(B,\mathbb{R}^2)$.
\end{lemma}

\begin{proof}
By Lemma 4 of \cite{vod-isan-MA} the kernel of the operator $Q$ on
the Sobolev class $W_{p,\loc}^1(\mathbb{H},\mathbb{C})$, $p> 1$, is
5-dimensional$:$ $u\in \kernel Q$ if and only if
\begin{align*}
u(z,t)&=a+ikz+tb+iz^2\overline{b}+i|z|^2b, \quad \text{where } a,b
\in \mathbb{C},\ k\in \mathbb{R}.
\end{align*}

Since the kernel of $Q$ is finite-dimensional the coercive estimate
of \cite{rom,rom1} is valid (see also \cite[Theorem 1]{isan-vod-EMJ}). 
Let $p>1$.
There are a constant $C=C(p)>0$ and a projection $P\colon W^1_{p}(B,\mathbb C)\to \ker Q$ such that
$$
\|D_hu-D_hP u\|_{p,B}\le C\|Qu\|_{p,B }
$$
for any $u\in W_p^1(B ,\Bbb{C})$.
 By standard way, we can use any projection operator in
the coercive estimate (see, for example, \cite[Proposition 2]{isan-smj2}).
Construct a projection $P$ on the kernel of $Q$. Let
$u_1,\dots,u_5\colon \mathbb{H} \to \mathbb{C}$ be the orthonormal
basis of $\operatorname{ker} (Q)$ with respect to the scalar product
$\langle \cdot,\cdot\rangle$ in the space
$L_2(B(\mathbf{0},1),\mathbb{C})$. Then set
$$
Pu=\sum_{i=1}^5 \langle u,u_i\rangle u_i.
$$
Obviously, $$
\|D_h Pu\|_{p,B}\leqslant C \sup\limits_{\mathbf{x}\in
B}|u(\mathbf{x})|.$$
The lemma follows.
\end{proof}

\subsection{Mappings with bounded specific oscillation}\label{SS:BSO}\label{subsection:BSO}

Here we give one result from \cite{isan-smj1} about mappings with bounded specific oscillation. 

\begin{lemma}[\mbox{\cite[Corollary to Theorem 1]{isan-smj1}}] \label{lem:BSO}
Let $U\subset\mathbb H$ be open, $q>1$. Suppose a mapping $f $ maps $U$ to the space of $2\times 2$-matrices and
for every ball $B\subset U$ there is orthogonal matrix $\phi_B\in SO(2)$ such that
\begin{equation}\label{eq:def_BSO}
\int\limits_B
|f(\mathbf x) - \phi_B|^q
d\mathbf x \leqslant \sigma^q \int\limits_B|\phi_B|^q d\mathbf x =\sigma^q |B |
\end{equation}
with some constant $\sigma>0$.

Then there is $\sigma_0>0$ s.t. 
$f \in L_{p,loc}(U)$ for all $p\in[q, 2\sigma_0/
\sigma)$ if $\sigma  < \sigma_0$. 
Moreover, if $q < p < \frac{(2-\delta)\sigma_0}{
\sigma}$, $\delta > 0$, then
$$
\int\limits_{B'} |f(\mathbf x) - \phi_{B}|^p d\mathbf x \leqslant C 
\sigma^{p-q}
\int\limits_{B'} |f(\mathbf x) - \phi_{B}|^q
d\mathbf x
$$
for every ball $B=B(\mathbf x,r)\subset U$, $B'=B\bigl(\mathbf x,\frac{9}{10}r\bigr)$.
Constant $C$ depends on $p,q, \delta$. 
\end{lemma}

Relation \eqref{eq:def_BSO} means that $f$ is a
mapping with bounded specific oscillation in $L_q$ with respect to $SO(2)$ ($f\in BSO_q(SO(2))$).
Mappings with bounded specific oscillation have common features with well-known $BMO$ class. Lemma \ref{lem:BSO} shows that, similar to $BMO$, mappings of $BSO_q(SO(2))$ class have some self-improving integrability property.

\subsection{Local quantitative Sobolev rigidity. Proof of Theorem \ref{th:local_Sob}}

In this subsection we will prove local quantitative Sobolev rigidity formulated in Theorem \ref{th:local_Sob}.

We have $F=(f_1,f_2,f_3)\in QI_{1+\varepsilon}(B(\mathbf 0,4))$. Denote $B=B(\mathbf 0,1)$.
Theorem \ref{th:local_d^H} implies
$$\sup\limits_{\mathbf x\in B}d^H(F(\mathbf x),\psi(\mathbf x))\leqslant N_3\varepsilon$$ for some isometry $\psi$. 
Without loss of generality we may assume that $\psi=\mathrm{id}$. Otherwise, we take $\psi^{-1}\circ F$ instead of $F$.

We also have 
qualitative rigidity in view of Lemma \ref{lem:Local_Sob_stab}:
\begin{equation*}
\sup\limits_{\mathbf x\in B} \rho( F(\mathbf x),\Phi(\mathbf x))
\leqslant \mu_1(\varepsilon,1/2),
\quad
\|D_h F-D_h\Phi\|_{2,B}
\leqslant
|B|^{1/2}\mu_2 (\varepsilon,1/2)
\end{equation*}
for some isometry $\Phi$. Here $\mu_i(\varepsilon,1/2)\to 0$ as $\varepsilon\to 0$, $i=1,2$.

Since $\sup\limits_{\mathbf x\in B} d^H(\mathbf x,\Phi(\mathbf x))\leqslant N_3\varepsilon+\mu_1(\varepsilon,1/2)$, Lemma \ref{lem:sB_isom} states
$$|D_h\Phi-I|\leqslant 2(N_3\varepsilon+\mu_1(\varepsilon,1/2))$$
and, therefore, 
$$\|D_h F-I\|_{2,B}
\leqslant
|B|^{1/2}\bigl(\mu_2 (\varepsilon,1/2)+2N_3\varepsilon+2\mu_1(\varepsilon,1/2)\bigr)
\overset{\mathrm{df}}{=}|B|^{1/2} \widetilde{\mu}_2(\varepsilon).$$

Denote  $B'=B(\mathbf 0,9/10)$, $u(\mathbf x)=(f_1(\mathbf x)-x_1,f_2(\mathbf x)-x_2)\in W^1_2(B',\mathbb R^2)$, $\mathbf x=(x_1,x_2,x_3)$. Obviously, $D_hu=D_hF-I$ and $|u(\mathbf x)|=d^H(\mathbf x,F(\mathbf x))$. Applying Lemma \ref{lem:coercive} to $u$, we obtain
\begin{equation}\label{eq:1}
\|D_h u\|_{2,B'}
\leqslant
C(\sup\limits_{\mathbf x\in B'}|u(\mathbf x)|+\|Qu\|_{2,B'})
\leqslant
C\varepsilon
+
C \|Q u\|_{2,B'}.
\end{equation}

Estimate $\|Qu\|_{2,B'}$. 
 The inequality
\begin{equation}\label{eq:main}
|Q (\mathbf{x}^{-1}\cdot F(\mathbf{x}))|\leqslant \frac{\varepsilon(\varepsilon+2)}{2}\bigl(|D_h
F(\mathbf{x})-I|+2\bigr) +\frac{1}{2}|D_h F(\mathbf{x})-I|^2
\end{equation}
holds almost everywhere in $B(\mathbf 0,4)$ \cite[Lemma 3]{vod-isan-MA} provided $F$ preserve $KR$-orientation. Lemma \ref{lem:QI->KR} guarantees that~$F$ preserves $KR$-orientation if $\varepsilon<\varepsilon_1$ and 
$$\sup\limits_{\mathbf x\in B}d^H(F(\mathbf x),\mathbf x)\leqslant N_3\varepsilon <1/4.$$

Take $\varepsilon<\min\{\varepsilon_1, (4N_3)^{-1}\}$. Then $F$ preserves $KR$-orientation and inequality \eqref{eq:main} is valid. Taking $L_2$-norm of both sides of \eqref{eq:main},
we obtain
\begin{equation}\label{eq:3}
\|Qu\|_{2,B'}\leqslant C\varepsilon+\frac{1}{2}\biggr(\int\limits_{B'} |D_hu(\mathbf x)|^4d\mathbf x\biggr)^2.
\end{equation}

It rests to estimate $\|D_h u\|^2_{4,B'}$. Lemma \ref{lem:Local_Sob_stab} allows us to apply Lemma \ref{lem:BSO}. Indeed, Lemma \ref{lem:Local_Sob_stab} states that for every ball $G\subset B(\mathbf 0,2)$ there is an isometry $\Psi_G$ such that
$$
\|D_h F-D_h\Psi_G\|_{2,G}\leqslant \mu_2(\varepsilon,1/2)|G|^{1/2}\leqslant \widetilde{\mu}_2(\varepsilon)|G|^{1/2}.
$$
Recall $\Psi_B=I$. Lemma  \ref{lem:BSO} implies
\begin{equation}\label{eq:2}
\|D_hu\|_{4,B'}^2=\|D_hF-I\|_{4,B'}^2\leqslant C \widetilde{\mu}_2(\varepsilon)\|D_hF-I\|_{2,B'}=C \widetilde{\mu}_2(\varepsilon)\|D_hu\|_{2,B'}
\end{equation}
if $\widetilde{\mu}_2(\varepsilon)<\sigma_0/3$.

Relations \eqref{eq:1}--\eqref{eq:2} imply
$$
\|D_hu\|_{2,B'}\leqslant C_1\varepsilon+C_2\widetilde{\mu}_2(\varepsilon)\|D_hu\|_{2,B'}.
$$
Finally,  $\|D_h u\|_{2,B'}
\leqslant 2C_1\varepsilon$ provided $\varepsilon$ is close enough to 0. It suffices to take $\varepsilon<\varepsilon_4$ where $\varepsilon_4<\min\{\varepsilon_1, (4N_3)^{-1}\}$, $C_2\widetilde{\mu_2}(\varepsilon_4)<1/2$, and $\widetilde{\mu_2}(\varepsilon_4)<\sigma_0/3$. 

If $\varepsilon\geqslant \varepsilon_4$ then the theorem is obvious:
$$
\|D_hF-I\|_{2,B'}\leqslant (2+\varepsilon)|B'|^{1/2}\leqslant \varepsilon \Bigl(\frac{2}{\varepsilon_4}+1\Bigr)|B'|^{1/2}.
$$

\subsection{Global quantitative Sobolev rigidity}

To pass from local to global we use the following proposition.

\begin{proposition}[\mbox{\cite{vod-isan-MA}}]\label{prop:Sob_norm:local_to_global}
Suppose $U$ is a John domain in $\mathbb{H}$ with inner radius
$\alpha$, outer radius $\beta$ and distinguished point $\mathbf x_*$, $\varkappa\geqslant 1$, $\sigma>0$. Let a
mapping $F\in QI_L( U)$ be such that for each ball
$B=B(\mathbf{a},r)$, $B(\mathbf{a},\varkappa r)\subset U$, there is
an isometry $\Phi_B$ meeting the condition
$$
\int_{B}
|D_h F(\mathbf{x})-D_h \Phi_B(\mathbf{x})|^2 d\mathbf{x} \leqslant \sigma^2 |r|^4.
$$
Then 
$$
\int_{U}\exp\Bigl(\Bigl(\frac{\beta}{\alpha}\Bigr)^{5}
\frac{c_2|D_hF(\mathbf{x})-D_h\Phi_0(\mathbf{x})|}{\sigma}\Bigr)\,d\mathbf{x}
\leqslant 16 |U|,
$$
where $\Phi_0=\Phi_{B_0}$, $B_0=B(\mathbf x_*,d_U(\mathbf x_*)/\varkappa)$.
Constant $c_2$ is independent of $U$ and $F$.
\end{proposition}

Proposition \ref{prop:Sob_norm:local_to_global}  shows that $BMO$ function $D_hF$ is exponentially integrable. The difference from well-known result is that we follow the coefficients of John domain.

The proof of Proposition \ref{prop:Sob_norm:local_to_global} follows word-by-word the part of the proof of Theorem 1 in \cite[Section 5]{vod-isan-MA}.

\end{document}